\documentclass[12pt,a4paper]{article}

\usepackage{algorithm}
\usepackage{algpseudocode}
\usepackage[utf8]{inputenc}
\usepackage[english]{babel}
\usepackage{color}
\usepackage{amsmath}
\usepackage{amsthm}
\usepackage{amsfonts}
\usepackage{amssymb}
\usepackage{graphicx}
\usepackage{mathtools}
\usepackage[left=2cm,right=2cm,top=2cm,bottom=2cm]{geometry}
\usepackage{hyperref}

\newcommand{\init}{\mathrm{Init}}
\newcommand{\unsafe}{\mathrm{Unsafe}}
\newcommand{\RR}{\mathbb{R}}
\newcommand{\ceni}{c_{\mathrm{I}}}
\newcommand{\cenu}{c_{\mathrm{U}}}
\newcommand{\elli}{E_{\mathrm{I}}}
\newcommand{\ellu}{E_{\mathrm{U}}}
\newcommand{\mpartial}[2]{\frac{\partial {#1}}{\partial {#2}}}
\newcommand{\dd}{\mathrm{d}}
\newcommand{\mder}[2]{\frac{\dd {#1}}{\dd {#2}}}

\newcommand{\LL}{\mathcal{L}(\chi, \lambda)}
\newcommand{\tbt}[4]{
\begin{bmatrix}
{#1} & {#2} \\
{#3} & {#4}
\end{bmatrix}
}
\newcommand{\tbts}[3]{
\begin{bmatrix}
{#1} & {#2} \\
{#2}^T & {#3}
\end{bmatrix}
}

\newtheorem{theorem}{Theorem}
\newtheorem{lemma}{Lemma}
\newtheorem*{problem}{Problem}

\title{
Solving Reachability Problems by a Scalable Constrained Optimization Method
 \thanks{This work was supported by the Czech Science Foundation (GACR) grant number 15-14484S with institutional support RVO:67985807.}}

\author{Jan Ku\v{r}\'{a}tko\thanks{ORCID: \href{http://orcid.org/0000-0003-0105-6527}{0000-0003-0105-6527}; Faculty of Mathematics and Physics, Charles University, Czech Republic; Institute of Computer Science, The Czech Academy of Sciences }, Stefan Ratschan\thanks{ORCID: \href{https://orcid.org/0000-0003-1710-1513}{0000-0003-1710-1513}; Institute of Computer Science, The Czech Academy of Sciences }}

\begin{document}
\maketitle
\abstract{In this paper we consider the problem of finding an evolution of a dynamical system that originates and terminates in given sets of states. However, if such an evolution exists then it is usually not unique. We investigate this problem and find a scalable approach for solving it. To this end we formulate an equality constrained nonlinear program that addresses the non-uniqueness of the solution of the original problem. In addition, the resulting saddle-point matrix is sparse. We exploit the structure in order to reach an efficient implementation of our method. In computational experiments we compare line search and trust-region methods as well as various updates for the Hessian.}
\section{Introduction}
\label{sec:Int}
In this paper we are concerned with the task of finding an evolution of a dynamical system that starts in a given set and reaches another set of states. We assume those two sets to be ellipsoids and call them the set of initial states and the set of unsafe states. Note that this problem does not have a unique solution, hence, it is not a classical boundary value problem (BVP). We solve it by formulating it as an equality constrained nonlinear programming problem. To this end we apply the \emph{sequential quadratic programming} method (SQP)~\cite{Nocedal:2006}. 

We will present an analysis of the resulting optimization problem and study the structure and properties of the Hessian matrix. In particular we are interested in the spectrum of the Hessian. We discuss which solution approaches are most suitable with respect to the properties of the Hessian and the whole saddle-point matrix. In addition we do computational experiments with both line-search SQP~\cite{Nocedal:2006} and trust-region SQP~\cite{Nocedal:2006} on benchmark problems.

The motivation for this work stems from the field of computer aided verification \cite{Branicky:2006,Lamiraux:2004,Zutshi:2013}. Here an evolution of a system from the set of initial states that reaches the set of unsafe states may represent a flaw in its design. To this end researchers try to develop automatized methods for identifying such flaws \cite{S-TaLiRo:2011}. There are several approaches to this problem~\cite{Abbas:2011,Zutshi:2013}.

The main contributions of this paper are:
\begin{itemize}
\item Formulation of the equality constrained nonlinear minimization problem that addresses the problem of nonunique solutions;
\item Investigation of the properties of the resulting saddle-point matrix such as:
\begin{itemize}
\item description of the structure of nonzero elements,
\item analytical formulas for the Hessian matrix and its use in the computation,
\item analysis of the spectrum of the Hessian for linear ODEs;
\end{itemize}
\item Comparison of line-search SQP with trust-region SQP.
\end{itemize}
The structure of the presented paper is the following. We formulate the problem we try to solve in Section~\ref{sec:PF}. In Section~\ref{sec:NPP} we formulate the nonlinear programming problem and review the SQP method in Section~\ref{sec:SPM}. In addition we describe the structure of nonzero elements in the saddle-point matrix and the properties regarding the spectrum of the Hessian in Section~\ref{sec:PofSPM}. The solution of the saddle-point system as well as the choice of the step length in the optimization process is left for Section~\ref{sec:StNLP}. Then there follow computational experiments in Section~\ref{sec:CE} and the paper finishes with a summary of results.
\section{Problem Formulation}
\label{sec:PF}
Consider a dynamical system whose dynamics is governed by a system of ordinary differential equations of the following form
\begin{equation}
	\label{eq:DiffEq}
	\dot{x} = f(x(t)),
\end{equation}
where $x: \RR^{\geq 0} \to \RR^n$ is an unknown function and the right hand side $f: \RR^n\to\RR^n$ is  continuously differentiable. We denote by $\Phi : \RR\times \RR^n \to \RR^n$ the \emph{Flow} function. Then for the initial time $t = 0$ we have $\Phi(0, x_0) = x_0$ and for $t \geq 0$, $\Phi(t, x_0) = x(t)$. We try to solve the following classical problem from the field of computer aided verification \cite{Branicky:2006, Lamiraux:2004, Zutshi:2013}.
\begin{problem}
Consider the dynamical system~\eqref{eq:DiffEq} and let $\init$ and $\unsafe$ be two sets of states in $\RR^n$. Find a solution of~\eqref{eq:DiffEq} and time $t \geq 0$ such that the initial state $x_0 \in \init$ and $\Phi(t, x_0) \in \unsafe$.
\end{problem}
\noindent
We assume that there exists such a solution and that the sets $\init$ and $\unsafe$ are disjoint. In addition, we assume the sets $\init$ and $\unsafe$ to be ellipsoids with centres $c_I \in \RR^n$ and $c_U \in \RR^n$, that is
\begin{align*}
\init & = \left\{ v \in \RR^n \mid \left( v - c_I\right)^T E_I( v- c_I) \leq 1\right\}, \\
\unsafe & = \left\{ v \in \RR^n \mid \left( v - c_U\right)^T E_U( v- c_U) \leq 1\right\}.
\end{align*}
We denote the norms induced by symmetric definite matrices $E_I \in \RR^{n \times n}$ and $E_U \in \RR^{n \times n}$ by $\| \cdot \|_{E_I}$, and by $\| \cdot \|_{E_U}$ respectively. Note that the problem we try to solve is a boundary value problem (BVP)  with separated boundary value conditions, however, it is not in standard form~\cite{Ascher:1998}:
\begin{itemize}
\item The upper bound $T$ on time $t \geq 0$ is unknown.
\item The boundary conditions are of the form
\[
g\left(x_0, t\right) =
\begin{bmatrix}
\| x_0 - c_I \|_{E_I} - 1 \\
\| \Phi(t, x_0) - c_U \|_{E_U} - 1
\end{bmatrix}
\leq
0,
\]
therefore, $g: \RR^{n+1} \to \RR^2$ yields an underdetermined BVP.
\end{itemize}
The unknown upper bound on time $T$ can be eliminated by transforming the BVP into an equivalent one with a fixed upper bound, introducing one more variable~\cite{Ascher:1995,Ascher:1998}. However, the problem we try to solve remains underdetermined.
%
%
%
\section{Nonlinear Programming Formulation}
\label{sec:NPP}
In order to find a solution of the dynamical system~\eqref{eq:DiffEq} from $\init$ to $\unsafe$ we reformulate it as a constrained minimization problem
\begin{equation}
\label{eq:minprob}
\min F(\chi) \qquad \text{subject to} \quad c_E(\chi) = 0 \text{ and } c_I(\chi) \leq 0,
\end{equation}
where $F:\RR^k \to \RR$, $c_E:\RR^k \to \RR^{m_E}$, $c_I:\RR^k \to \RR^{m_I}$ are twice continuously differentiable.

\begin{figure}
\centering
\includegraphics[scale=0.8]{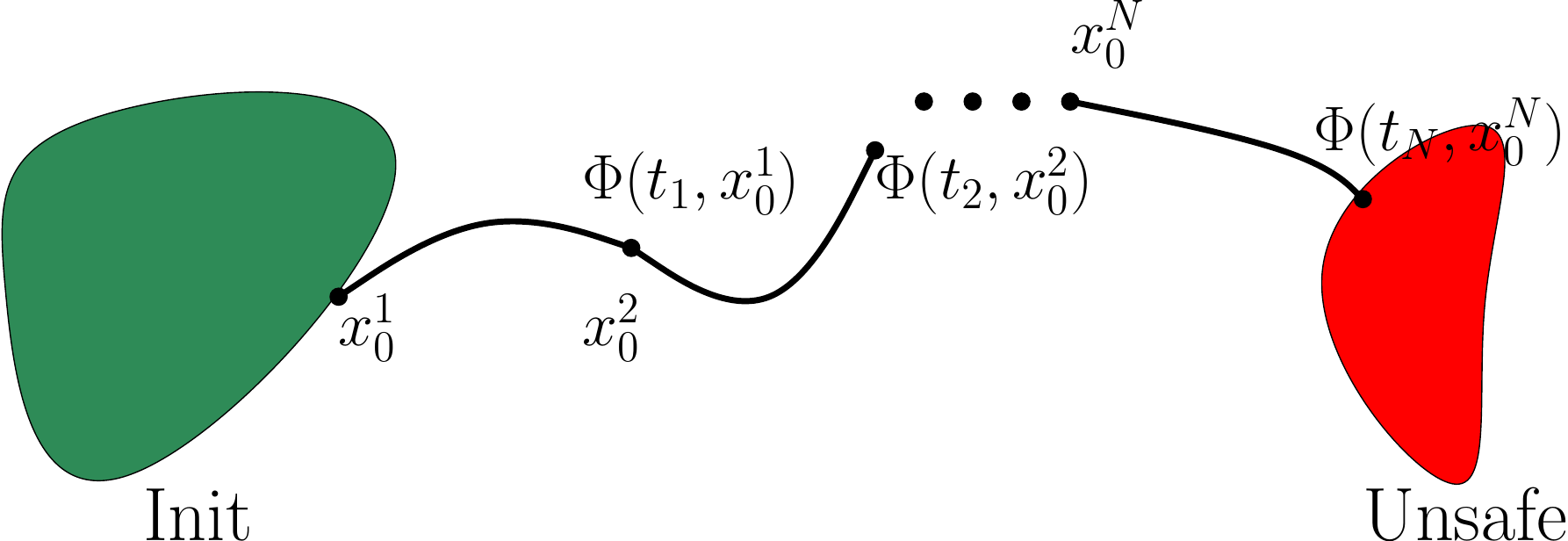}
\caption{Illustration of connected $N$ solution segments forming one solution of the dynamical system~\eqref{eq:DiffEq} from $\init$ to $\unsafe$.}
\label{fig:solution_ilustrated}
\end{figure}

To solve the minimization problem~\eqref{eq:minprob} we use the idea of multiple shooting~\cite{Ascher:1995} where one finds a solution by patching up a number of solution segments. One such a solution to our problem stated above is illustrated in~Fig. \ref{fig:solution_ilustrated}. The unknown vector $\chi$ consists of the initial states of solution segments $x_0^i \in \RR^n$ and their lengths $t_i \geq 0$, $1 \leq i \leq N$, and has the following form
\begin{equation}
\label{eq:parameters}
\chi = \left[ x_0^1, t_1, x_0^2, t_2, \ldots, x_0^N, t_N \right]^T \in \RR^{N(n+1)}\ .
\end{equation}
The problem has infinitely many solutions since our only requirement on the lengths $t_i$, $1 \leq i \leq N$, is that their sum $\sum_i t_i = t$. Therefore, if one sets $x_0^1 \in  \init$ and $\Phi(t_N, x_0^N) \in \unsafe$ as inequality constraints $c_I(\chi) \leq 0$ and matching conditions $g_i(\chi) := x_0^{i+1} - \Phi(t_i, x_0^i) = 0$, $1 \leq i \leq N-1$, as equality constraints $c_E(\chi)$, then one still is in need of regularization to address the problem of having infinitely many solutions.

We studied various formulations of the objective function $F(\chi)$ and investigated the possibility of adding a regularization term. Our goal was to place additional conditions on the lengths $t_i$, $1 \leq i \leq N$, of solution segments. To this end we use the objective function
\begin{equation}
\label{eq:OF}
F( \chi) = \frac{1}{2}\sum_{i = 1}^N t_i^2.
\end{equation}
This choice of the objective function~\eqref{eq:OF} has interesting consequences and the idea behind them is shown in the following Lemma~\ref{lem:shortest_time_idea}.
\begin{lemma}
\label{lem:shortest_time_idea}
Let $L$ be a positive scalar and $N \in \mathbb{N}$. Then the solution to
\[
\min_{t} \frac{1}{2}\sum_{i=1}^N t_i^2\quad \text{subject to}\quad L = \sum_{i=1}^Nt_i
\]
is unique and attained at $t_1 = t_2 = \cdots = t_N = L/N$.
\begin{proof}
The uniqueness of the solution follows from the convexity of the objective function and the linearity of the constraint. When one forms the Lagrangian $\mathcal{L}(t,u)$ and differentiate with respect to $t$, then $\nabla_t \mathcal{L}(t,u) = [t_1-u, \ldots, t_N - u]^T \in \RR^{N}$. Hence, putting $\nabla_t \mathcal{L}(t,u) = 0$ gives $t_1 = t_2 = \cdots = t_N = u$, where $u \in \RR$ is the Lagrange multiplier. The rest follows from the fact that $t_i$, $1 \leq i  \leq N$, are equal and their sum is equal to $L$.
\end{proof}
\end{lemma} 
We observed in our numerical experiments that the objective function~\eqref{eq:OF} drives the solution segments to be of equal length at the end of optimization process. In addition, the use of the objective function $F(\chi)$ resembles the shortest time problem~\cite{Betts:2010,Pontryagin:1962}. This leads to solutions of the minimization problem~\eqref{eq:minprob} for which $x_0^i \in \init$ and $\Phi(t_N, x_0^N) \in \unsafe$ are on the boundaries of $\init$, and $\unsafe$ respectively. 

Because of this phenomenon we formulate boundary conditions $x_0^1 \in \init$ and $\Phi(t_N, x_0^N) \in \unsafe$ as equality constraints. Furthermore, the objective function~\eqref{eq:OF} and constraints~\eqref{eq:constraints}
\begin{equation}
\label{eq:constraints}
c_E(\chi)  =
\begin{bmatrix}
x_0^2 - \Phi(t_1, x_0^1) \\
x_0^3 - \Phi(t_2, x_0^2) \\
\vdots \\
x_0^N - \Phi(t_{N-1}, x_0^{N-1}) \\
\end{bmatrix},
\quad
c_I(\chi)
=
\frac{1}{2}
\begin{bmatrix}
\|x_0^1 - \ceni \|_{\elli}^2 - 1 \\
 \|\Phi(t_N, x_0^N) - \cenu \|_{\ellu}^2 - 1  \\
\end{bmatrix}
\end{equation}
are separable. This will allow us to design efficient methods~\cite{Griewank:1982} for solving the optimization problem~\eqref{eq:minprob}. In particular, we exploit the separability in the computation and the approximation of the Hessian, see Theorems~\ref{lem:nonlinH} and~\ref{lem:linH}. Moreover, we use the separability in the implementation of the projected preconditioned conjugate gradient method as described in Section~\ref{subsec:FurRem}.

Before we chose the objective function~\eqref{eq:OF} we had considered and tried various  further formulations of the minimization problem~\eqref{eq:minprob}. In more detail we investigated
\begin{align}
\label{eq:MinGaps}
& \min \frac{1}{2} \sum_{i=1}^{N-1} \| x_0^{i+1} - \Phi(t_i, x_0^i)\|^2 \text{ subject to }   c_I(\chi) \leq 0, \\
\label{eq:MinDist}
& \min \frac{1}{2}\left( \|x_0^1 - c_I\|_{E_I}^2 + \|\Phi(t_N, x_0^N) - c_U\|_{E_U}^2\right) \text{ subject to } c_E(\chi) = 0, \\
\label{eq:MinAll}
& \min \frac{1}{2}\left( \|x_0^1 - c_I\|_{E_I}^2  + \|\Phi(t_N, x_0^N) - c_U\|_{E_U}^2 +\sum_{i=1}^{N-1}\| x_0^{i+1} - \Phi(t_i, x_0^i)\|^2 \right). 
\end{align}
However,  in formulations~\eqref{eq:MinGaps}-\eqref{eq:MinAll} we have no control over the lengths $t_i$, $1 \leq i \leq N$, of solution segments. This causes problems during the numerical computation where the length of one solution segment may be large and another solution segment may degenerate to zero length. Because of this very reason we prefer the objective function~\eqref{eq:OF}.

Note that adding the regularization term~\eqref{eq:OF} to~\eqref{eq:MinGaps}-\eqref{eq:MinAll} does not necessarily help. If one does that for the minimization problems~\eqref{eq:MinGaps} and~\eqref{eq:MinAll}, then the computed solution is not necessarily continuous since at the minimum $\sum_{i=1}^{N-1}\| x_0^{i+1} - \Phi(t_i, x_0^i)\|^2 \neq 0$.

There are other possibilities how to choose a regularization term such as $\sum (t_i - t_{i+1})^2$, and $\sum(\bar{t} - t_i)^2$ where $\bar{t}$ is the arithmetic mean of the lengths $t_i$, $1 \leq i  \leq N$. These also force lengths of solution segments to be equally distributed. However, the objective function and constraints are no longer as separable as in the case of using the objective function~\eqref{eq:OF}. 
%
%
%
\section{Review of SQP}
\label{sec:SPM}
In this section we review the SQP method and introduce the notation. In order to find a solution from $\init$ to $\unsafe$, we solve the following nonlinear minimization problem
\begin{equation}
\label{eq:MinProbEQ}
\min_\chi \frac{1}{2}\sum_{i=1}^Nt_i^2\quad \text{subject to}\quad c(\chi) = 0,
\end{equation}
where the vector of constraints $c(\chi)$ has the form
\begin{equation}
\label{eq:constrEQ}
c(\chi)
=
\begin{bmatrix}
\frac{1}{2}\left(  \|x_0^1 - \ceni \|_{\elli}^2- 1\right) \\
c_E(\chi) \\
\frac{1}{2}\left(  \| \Phi(t_N, x_0^N) - \cenu \|_{\ellu}^2- 1\right)
\end{bmatrix} \in \RR^{(N-1)n + 2}.
\end{equation}
When one solves \eqref{eq:MinProbEQ} by the SQP method~\cite{Nocedal:2006} the Lagrangian $\LL$ is formed, where 
$\lambda \in \RR^{(N-1)n+2}$ is a vector of Lagrange multipliers such that
\begin{equation}
\label{eq:lagmul}
\lambda = [\lambda_\mathrm{I}, \lambda_1, \ldots, \lambda_{N-1}, \lambda_\mathrm{U}]^T \in \RR^{(N-1)n+2}\,
\end{equation}
with $\lambda_\mathrm{I} \in \RR$, $\lambda_\mathrm{U} \in \RR$ and $\lambda_i \in \RR^{n}$ for $1 \leq i \leq N-1.$ Let us denote by $B(\chi) \in \RR^{N(n+1)\times (N-1)n+2}$ the Jacobian of the vector of constraints $c(\chi)$ and assume it has full column rank.  Then for the Lagrangian $\LL = F(\chi) + \lambda^Tc(\chi)$ one gets
\begin{align}
\label{eq:nondx}
\nabla_\chi \LL  & = \nabla_\chi F(\chi) + B(\chi)\lambda\,, \\
\label{eq:nondlam}
\nabla_\lambda \LL & = c(\chi)\,.
\end{align}
The solution $\chi^\star$ of the minimization problem~\eqref{eq:MinProbEQ} satisfies the Karush-Kuhn-Tucker (KKT) conditions if there exists $\lambda^\star$ such that $\nabla_\chi \mathcal{L}(\chi^\star, \lambda^\star) = 0$ and $\nabla_\lambda \mathcal{L}(\chi^\star, \lambda^\star) = 0$. We use the iterative method to solve the problem~\eqref{eq:MinProbEQ}. Then
the resulting saddle point system we need to solve in every iteration \cite{LuksanVlcek:2001,Nocedal:2006} is
\begin{align}
\label{eq:KKTsystem}
\tbts{H}{B}{0}
\begin{bmatrix}
d_\chi \\
d_\lambda
\end{bmatrix}
& =
-
\begin{bmatrix}
\nabla_\chi \LL \\
\nabla_\lambda \LL
\end{bmatrix},
\end{align}
which we write shortly as
\[
\notag
Kd = 
b\,,
\]
where by $H \in \RR^{N(n+1)\times N(n+1)}$ we denote $\nabla_\chi^2 \LL$. The solution vector of the saddle point system~\eqref{eq:KKTsystem} is then used to compute the next iterate 
\begin{equation}
\label{eq:nextit}
\chi^+ = \chi + \alpha_\chi d_\chi \textrm{ and } \lambda^+ = \lambda + \alpha_\lambda d_\lambda\,, 
\end{equation}
where the values $\alpha_\chi = \alpha_\lambda = 1$ give the Newton method. One can compute the initial value of $\lambda$~\cite[Ch.~18]{Nocedal:2006}, however, we set the initial value to be $\lambda = [1, \ldots, 1]^T \in \RR^{(N-1)n+2}$ for numerical experiments in this paper .
\section{Properties of the Saddle-point Matrix}
\label{sec:PofSPM}
In this section we discuss and show the structure of nonzero elements in the saddle-point matrix $K$. We show in Lemma~\ref{lem:Jac_LIcolumns} the linear independence of columns of the Jacobian $B$. In addition we show that the Hessian matrix $H$ is block-diagonal, indefinite and singular.
\subsection{Jacobian of Constraints}
\label{subsec:JoC}
\begin{figure}
\centering
\includegraphics[width = 0.6\textwidth]{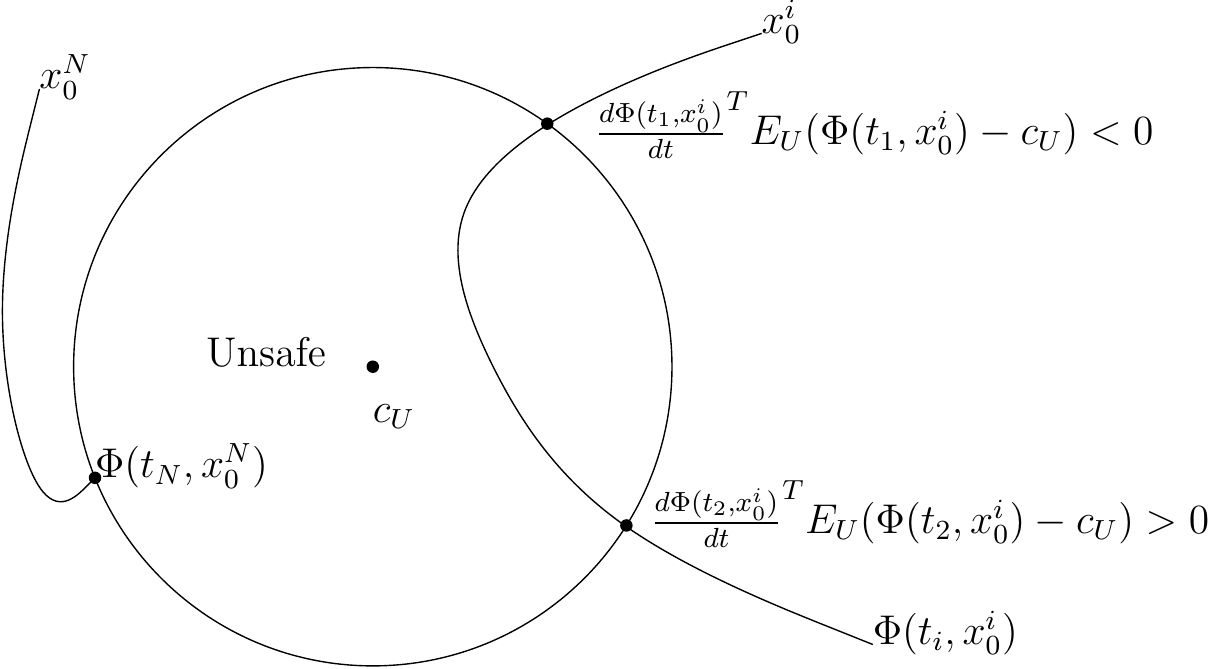}
\caption{Solution segments and their intersection with the boundary of the set $\unsafe$. At the point of entry to $\unsafe$ is  $\frac{d \Phi(t_N, x_0^N)}{d t_N}^T E_U  \left( \Phi(t_N, x_0^N) - c_U \right) < 0$, since it measures how the value $\frac{1}{2}\left( \| \Phi(t_N, x_0^N) - \cenu\|_{\ellu}^2 - 1\right)$ changes with respect to $t_N$ for the fixed initial state $x_0^N$.}
\label{fig:point_of_entry}
\end{figure}

First, we describe the Jacobian of constraints $B$. We denote the sensitivity function $S: \RR\times\RR^n \to \RR^{n \times n}$ of the solution $x(t)= \left[ x_1(t), \ldots, x_n(t) \right]^T \in \RR^n$ of the differential equation~\eqref{eq:DiffEq} to the change of the initial value $x_0 = \left[x_{0,1}, \ldots, x_{0,n} \right]^T \in \RR^n$ by
\begin{equation}
\label{eq:Sensitivity}
S(t, x_0) =  
\begin{bmatrix}
 \frac{\partial x_1(t)}{\partial x_{0,1}} & \ldots & \frac{\partial x_1(t)}{\partial x_{0,n}} \\
\vdots & \ddots & \vdots \\
 \frac{\partial x_n(t)}{\partial x_{0,1}} & \ldots & \frac{\partial x_n(t)}{\partial x_{0,n}} \\
\end{bmatrix}.
\end{equation}
%
\begin{lemma}
\label{lem:BIBE_Jacobians}
The Jacobians of the vectors of constraints~\eqref{eq:constraints} are
\[
B_I(\chi)
=
\begin{bmatrix}
E_I (x_0^1 - c_I)& 0 \\
0 & \vdots \\
\vdots & S(t_N, x_0^N)^T E_U \left( \Phi(t_N, x_0^N) - c_U \right)  \\
0 & \frac{d \Phi(t_N, x_0^N)}{d t_N}^T E_U  \left( \Phi(t_N, x_0^N) - c_U \right)
\end{bmatrix} \in \RR^{N(n+1) \times 2}
\]
and
\[
B_E(\chi) 
= 
\begin{bmatrix}
 -S(t_1,x_0^1)^T  & &   & &  \\
-\frac{d \Phi(t_1, x_0^1)}{d t_1}^T &   &&  &  \\
I &  - S(t_2, x_0^2)^T & &  & \\
  & -\frac{d \Phi(t_2, x_0^2)}{d t_2}^T & & &  \\
  &  I &  & &  \\
 &   & \ddots&  & \\
  & & & &  \\
 &   &  &   & -S(t_{N-1},x_0^{N-1})^T   \\
 &  & &  &  -\frac{d \Phi(t_{N-1}, x_0^{N-1})}{d t_{N-1}}^T  \\
 &  & &  & I\\ 
  &  & &  & 0
\end{bmatrix} \in \RR^{N(n+1)\times (N-1)n},
\]
where $I \in \RR^{n \times n}$ is the identity matrix. Matrices $S(t_i,x_0^i) \in \RR^{n \times n}$ are sensitivities and $\Phi(t_i, x_0^i)$ are final states of the $i$th solution segment of~\eqref{eq:DiffEq} with the initial state $x_0^i$ and the length $t_i > 0$, $1 \leq i \leq N$. 
\begin{proof}
Let us start with $B_I(\chi)$. The first constraint $g_I(\chi) = \frac{1}{2}\left((x_0^1 - c_I)^TE_I(x_0^1 - c_I)-1\right)$ depends only on $x_0^1 \in \RR^n$, therefore $\nabla_\chi g_I(\chi) = [(x_0^1-c_I)^TE_I, 0, \ldots, 0]^T \in \RR^{N(n+1)}$. The second constraint $g_U(\chi) = \frac{1}{2}\left((\Phi(t_N, x_0^N) - c_U)^TE_U(\Phi(t_N, x_0^N) - c_U)-1\right)$ depends both on $x_0^N \in \RR^n$ and $t_N \in \RR$. Using the chain rule one obtains after differentiating with respect to the initial state $\partial g_I(\chi)/\partial x_0^N = S(t_N, x_0^N)^TE_U(\Phi(t_N, x_0^N) - c_U) \in \RR^n$ and with respect to time $d g_U(\chi)/ d t_N =\left( d \Phi(t_N, x_0^N)/ d t_N \right)^T E_U  \left( \Phi(t_N, x_0^N) - c_U \right) \in \RR$.

To get $B_E(\chi)$ one needs to differentiate $g_i(\chi)$, $1 \leq i \leq N-1$. Because of the ordering of parameters in~\eqref{eq:parameters} one obtains the banded structure of the Jacobian $B_E(\chi)$. Constraints $g_i(\chi) = x_0^{i+1} - \Phi(t_i, x_0^i) \in \RR^n$ depend on $x_0^{i+1} \in \RR^n$, $x_0^i \in \RR^n$ and $t_i \in \RR$ for $1 \leq i \leq N-1$. Therefore $\partial g_i(\chi)/\partial x_0^{i+1} = I \in \RR^{n \times n}$ is the identity matrix, $d g_i(\chi)/d t_i = - (d\Phi(t_i, x_0^i)/d t_i)^T$ and $\partial g_i(\chi)/\partial x_0^i = - S(t_i, x_0^i)^T$. No constraint in $c_E(\chi)$ depends on $t_N$, therefore the last row of $B_E(\chi)$ is the zero vector.
\end{proof}
\end{lemma}
 When $E_I(x_0^1 - c_I) \in \RR^n$ is nonzero and $ \frac{d \Phi(t_N, x_0^N)}{d t_N}^T E_U  \left( \Phi(t_N, x_0^N) - c_U \right) \neq 0$, then the Jacobian $B$ has full column rank, as shown in Lemma~\ref{lem:Jac_LIcolumns}. Assuming the term $E_I(x_0^1 - c_I)$ to be nonzero is natural in the sense that the choice of the objective function~\eqref{eq:OF} and constraints~\eqref{eq:constrEQ} place $x_0^1$ on the boundary of $\init$, hence $x_0^1 \neq c_I$. Similarly, the final state of the last solution segment $\Phi(t_N, x_0^N)$ is placed on the boundary of $\unsafe$. The value of $\frac{d \Phi(t_N, x_0^N)}{d t_N}^T E_U  \left( \Phi(t_N, x_0^N) - c_U \right) \leq 0$ at the point of entry to the set $\unsafe$ is illustrated in Fig.~\ref{fig:point_of_entry}. 
 \begin{lemma}
\label{lem:Jac_LIcolumns}
Let $A \in \RR^{N(n+1)\times ((N-1)n+2)}$ be a matrix of the form
\[
A = \begin{bmatrix}
w_1 & M_1 & &   & & & \\
&v_1^T &   &&  & & \\
&I &  M_2 & & & & \\
 & & v_2^T & & & & \\
 & &  I & \ddots & & & \\
& &   & \ddots & v_{N-2}^T &  & \\
& &  & & I & M_{N-1}  & \\
& &  & &  & v_{N-1}^T & \\
& &  & &  & I & w_2\\ 
 & &  & &  & 0 & \omega\\ 
\end{bmatrix}\ ,
\]
where for $1 \leq i \leq N-1$, $M_i \in \RR^{n \times n}$, vectors $v_i, w_1, w_2 \in \RR^n$, $\omega \in \RR$, and $I \in \RR^{n \times n}$ is the identity matrix. If $w_1 \in \RR^n$ is a non-zero vector and $\omega \neq 0$, then the matrix $A$ has full-column rank.
\begin{proof}
We prove this Lemma by contradiction. Suppose columns in $A$ are linearly dependent, therefore, there exists a non-zero vector $y = \left[ \alpha, y_1^T, \ldots, y_{N-1}^T, \beta \right]^T \in \RR^{(N-1)n+2}$ with $y_i \in \RR^n$, $1 \leq i \leq N-1$, $\alpha \in \RR$, and $\beta \in \RR$ so that
\begin{align*}
\alpha w_1 + M_1y_1 & = 0,\\
y_i + M_{i+1}y_{i+1} & = 0,   \\
y_{N-1} + \beta w_2 & = 0, \\ 
\beta\omega = 0 ,
\end{align*}
for $1 \leq i \leq N-2$, and $v_i^Ty_i  = 0$ for $1 \leq i \leq N-1$. Since we assume $\omega$ to be a non-zero scalar, therefore, we get $\beta = 0$. It follows that $y_{N-1} = 0 \in \RR^n$. If we substitute into formulae above we obtain $y_i = 0 \in \RR^n$ for $1 \leq i \leq N-2$. Therefore, also $\alpha w_1 = 0 \in \RR^{n}$. For $\alpha \neq 0$ this is only possible if $w_1 = 0 \in \RR^{n}$. This is contradiction with the assumption that $w_1$ is a non-zero vector. 
\end{proof}
\end{lemma}
\subsection{Hessian Matrix}
\label{subsec:HM}
The following theorems describe the block-diagonal structure of the Hessian $\nabla_\chi^2\LL$ depending on the dynamics of \eqref{eq:DiffEq}.
\begin{theorem}
\label{lem:nonlinH}
Consider the dynamical system \eqref{eq:DiffEq} and the minimization problem \eqref{eq:MinProbEQ}. Then the Hessian of the Lagrangian is  block diagonal of the form
\[
\nabla_\chi^2\LL
=
\begin{bmatrix}
\tbts{A_1}{v_1}{\alpha_1} & & \\
& \ddots & \\
& & \tbts{A_N}{v_N}{ \alpha_N}
\end{bmatrix} \in \RR^{N(n+1) \times N(n+1)}\,,
\]
where $n$ is the statespace dimension, and $N$ is the number of segments. For $1 \leq i \leq N$, blocks $A_i \in \RR^{n \times n}$, $v_i \in \RR^n$ and $\alpha_i \in \RR$. Here
\begin{align*}
v_i & = -\left(\mpartial{f(x(t_i))}{x}S(t_i, x_0^i)\right)^T\lambda_{i}\,,\quad 1\leq i \leq N-1, \\
\alpha_i & = -\left(\mpartial{f(x(t_i))}{x} f(x(t_i))\right)^T\lambda_{i} + 1\,,\quad  1\leq i \leq N-1, \\
A_1 & = \lambda_{I}\elli - \frac{\partial^2 \Phi(t_1, x_0^1)}{\partial {x_0^1}^2} \circ \lambda_{1} \,, \\
A_i & = -\frac{\partial^2 \Phi(t_i, x_0^i)}{\partial {x_0^i}^2} \circ \lambda_{i}\,,\quad 2\leq i \leq N-1,
\end{align*}
where $\frac{\partial^2 \Phi(t_i, x_0^i)}{\partial {x_0^i}^2}  \in \RR^{n \times n \times n}$ is a tensor and the symbol $\circ$ denotes the contraction by $\lambda_{i}$. The last block of $\nabla_\chi^2 \LL$ consists of
\begin{align*}
A_N & = \lambda_{U}S(t_N, x_0^N)^T \ellu S(t_N, x_0^N) + \lambda_{U}\frac{\partial^2 \Phi(t_N, x_0^N)}{\partial {x_0^N}^2}^T\circ \ellu \left( \Phi(t_N, x_0^N) - \cenu \right)\,, \\
v_N & =\lambda_{U} \left( \mpartial{f(x(t_N))}{x}S(t_N, x_0^N) \right)^T \ellu \left( \Phi(t_N, x_0^N) - \cenu \right) +\lambda_{U}S(t_N, x_0^N)^T \ellu f(x(t_N))\,, \\
\alpha_N & =\lambda_{U} \left( \mpartial{f(x(t_N))}{x}f(x(t_N)) \right)^T \ellu \left( \Phi(t_N, x_0^N) - \cenu \right) +\lambda_{U}f(x(t_N))^T \ellu f(x(t_N)) + 1\,,
\end{align*}
where $\frac{\partial^2 \Phi(t_N, x_0^N)}{\partial {x_0^N}^2} \in \RR^{n \times n \times n}$ is a tensor.
\begin{proof}
The proof follows from differentiating $\LL$ twice with respect to the parameter $\chi$. First, one gets $\nabla_\chi \LL = \nabla_\chi F(\chi) + B(\chi)\lambda$, and then
\begin{align*}
\nabla_\chi^2 \LL & =   \nabla_\chi^2F(\chi) + \lambda_{I}\nabla_\chi^2g_\mathrm{I}(\chi) +  \\
& \quad  +  \lambda_{1} \circ \nabla_\chi^2g_1(\chi)  + \cdots + \lambda_{N-1} \circ \nabla_\chi^2g_\mathrm{N-1}(\chi) +  \lambda_{U}\nabla_\chi^2g_\mathrm{U}(\chi)\,.
\end{align*}
The term  $\nabla_\chi^2F(\chi) \in \RR^{N(n+1) \times N(n+1)}$ is a diagonal matrix containing the second derivatives with respect to $t_i$, $1 \leq i \leq N$, of the term $\frac{1}{2}\sum_i^N t_i^2$. Therefore, there are only $N$ nonzero elements $ \dd^2 F(\chi)/\dd t_i^2$, $1 \leq i \leq N$. Those are the ``$+1$'' in formulas for $\alpha_i$, $1 \leq i \leq N$.

Let us have a look at the matching conditions $g_i(\chi) = x_0^{i+1} - \Phi(t_i, x_0^i)$, $1 \leq i \leq N-1$, and observe that those are dependent on $x_0^{i+1}$, $x_0^i$ and $t_i$. The term $g_i(\chi)$ vanishes when differentiated twice with respect to $x_0^{i+1}$. Therefore $\lambda_{i} \circ \nabla_\chi^2g_i(\chi) \in \RR^{N(n+1) \times N(n+1)}$ features possibly nonzero entries only in elements corresponding to the second and mixed derivatives of $x_0^i$ and $t_i$, $1 \leq i \leq N-1$. Hence one obtains the block-diagonal structure of the Hessian.

To compute the second mixed derivatives first put $\dd \Phi(t_i, x_0^i)/\dd t = f(x(t_i)) \in \RR^n$ and then $\partial f(x(t_i))/\partial x_0^i = \partial f(x(t_i))/\partial x \cdot S(t_i, x_0^i) \in \RR^{n \times n}$. Those are contained in formulas for $v_i$, $1 \leq i \leq N-1$. The second derivative with respect to $t$, using the chain rule, is then $\dd^2 \Phi(t_i, x_0^i)/\dd t^2 = \dd f(x(t_i))/ \dd t = \partial f(x(t_i))/\partial x \cdot f(x(t_i)) \in \RR^n$. These formulas you can find in $\alpha_i$, $1 \leq i \leq N-1$.

In the case of the first constraint $g_\mathrm{I}(\chi) = 0.5\left((x_0^1 - \ceni)^T\elli(x_0^1 - \ceni)-1\right)$, its second derivative with respect to $x_0^1$ is $\elli \in \RR^{n \times n}$. The computation is more difficult for $g_\mathrm{U}(\chi) = 0.5\left((\Phi(t_N, x_0^N) - \cenu)^T\ellu (\Phi(t_N, x_0^N) - \cenu) - 1\right)$, however, $A_N$, $v_N$ and $\alpha_N$ is again the result of the application of the chain rule. The first derivatives with respect to $x_0^N$ and $t_N$ are 
\begin{align*}
\mpartial{g_\mathrm{U}(\chi)}{x_0^N} & = S(t_N, x_0^N)^T\ellu(\Phi(t_N, x_0^N)-\cenu)\,, \\ 
\mder{g_\mathrm{U}(\chi)}{t_N} & = f(x(t_N))^T\ellu (\Phi(t_N, x_0^N) - \cenu)\,.
\end{align*}
Differentiating both terms again with respect to $x_0^N$ and $t_N$ delivers desired formulas for  $A_N$, $v_N$ and $\alpha_N$.
\end{proof}
\end{theorem}
The \emph{Flow} of a linear ODE $\dot{x} = Ax(t)$ is $\Phi(t, x_0) = e^{At}x_0$, where $x_0$ is an initial state and $t \geq 0$. Then $S(t, x_0) = e^{At}$, the second derivative with respect to $x_0$ vanishes, and $f(x(t)) = Ae^{At}x_0$. From this one obtains the following theorem.
\begin{theorem}
\label{lem:linH}
Consider a linear dynamical system $\dot{x} = Ax(t)$ and the minimization problem \eqref{eq:MinProbEQ}. Then the Hessian of the Lagrangian is  block diagonal of the form
\[
\nabla_\chi^2\LL
=
\begin{bmatrix}
\tbts{A_1}{v_1}{\alpha_1} & & \\
& \ddots & \\
& & \tbts{A_N}{v_N}{ \alpha_N}
\end{bmatrix} \in \RR^{N(n+1) \times N(n+1)}\,,
\]
where $n$ is the statespace dimension, and $N$ is the number of segments. For $1 \leq i \leq N$, blocks $A_i \in \RR^{n \times n}$, $v_i \in \RR^n$ and $\alpha_i \in \RR$. Here
\begin{align*}
v_i & = -\left(Ae^{At_i}\right)^T\lambda_{i} \,,\quad 1\leq i \leq N-1, \\
\alpha_i & = -\left(A^2e^{At_i}x_0^i\right)^T\lambda_{i} + 1\,,\quad  1\leq i \leq N-1, \\
A_1 & = \lambda_{I}\elli\,, \\
A_i & =0\,,\quad 2\leq i \leq N-1, \\
A_N & = \lambda_{U} e^{A^T t_N} \ellu e^{At_N}\,, \\
v_N & =  \lambda_{U} \left( A e^{At_N}\right)^T \ellu \left( e^{At_N}x_0^N - \cenu\right) + \lambda_{U} e^{A^T t_N}\ellu Ae^{At_N}x_0^N \,, \\
\alpha_N & = \lambda_{U}\left( A^2 e^{At_N}x_0^N\right)^T \ellu \left( e^{At_N}x_0^N - \cenu\right) + \lambda_{U} \left( Ae^{At_N}x_0^N\right)^T \ellu Ae^{At_N}x_0^N + 1\,.
\end{align*}
\begin{proof}
The result follows directly from Theorem~\ref{lem:nonlinH}.
\end{proof}
\end{theorem}
Theorems~\ref{lem:nonlinH} and~\ref{lem:linH} show the block diagonal structure of the Hessian matrix $H$. In addition, when the dynamics in~\eqref{eq:DiffEq} is linear, one does not need to differentiate the Lagrangian $\LL$ twice with respect to $\chi$. Computed data from the Jacobian $B$, in particular sensitivity functions $S(t_i, x_0^i) = e^{At_i}$, $1 \leq i \leq N-1$, can be used in the formulas of Theorem~\ref{lem:linH}. Because of this observation, for linear systems~\eqref{eq:DiffEq}, one can work with the Hessian given by analytical formulas with no extra computational effort.

We are interested in the spectrum of $H$ and the whole saddle point matrix $K$ as well as the conditions on solvability of the saddle-point system~\eqref{eq:KKTsystem}. Since there is a complete description of $H$ in Theorem~\ref{lem:linH}, let us start there.
\begin{lemma}
\label{lem:blockHlin_eig}
Let $M \in \RR^{(n+1)\times (n+1)}$ be a matrix of the form
\[
M
=
\tbts{0}{v}{\alpha}\,,
\]
where the upper-left block $0 \in \RR^{n \times n}$, $v \in \RR^n$  and $\alpha \in \RR$. Then $M$ has $(n-1)$ zero eigenvalues and additional two such that
\[
\lambda_\pm = \frac{\alpha \pm \sqrt{\alpha^2 + 4v^Tv}}{2}\,.
\]
\begin{proof}
Eigenvalues and eigenvectors $(\lambda, u)$ of $M$ satisfy $Mu = \lambda u$ with $u = [x^T, y^T]^T$ where
\[
\tbts{0}{v}{\alpha}
\begin{bmatrix}
x \\
y
\end{bmatrix}
=
\lambda
\begin{bmatrix}
x \\
y
\end{bmatrix}\,.
\]
First assume that $\lambda = 0$. Then one can find $(n-1)$ orthogonal eigenvectors $u = [x^T, 0]^T$ satisfying $v^Tx = 0$ since $v \in \RR^n$. 

By rewriting the matrix equation above one gets
\begin{align*}
vy & = \lambda x \\
v^Tx + \alpha y & = \lambda y\,.
\end{align*}
When $\lambda \neq 0$ then  $y$ needs to be nonzero. It follows from the first equation that for $y = 0$ one gets $x = 0$. Dividing the first equation by $\lambda$ and substituting for $x$ in the second, one obtains $v^Tvy + \lambda \alpha y = \lambda^2 y$. Since $y \neq 0$ we can divide both sides and put $\lambda^2 - \lambda \alpha - v^Tv = 0$.
\end{proof}
\end{lemma}
\begin{lemma}
\label{lem:H1HN_lin}
Let $M \in \RR^{(n+1)\times(n+1)}$ be a matrix of the form
\[
M
=
\tbts{A}{v}{\alpha}\,,
\]
where $A \in \RR^{n \times n}$ is SPD (SND), $v \in \RR^n$ and $\alpha \in \RR$. Then $M$ has $n$ strictly positive (negative, respectively) eigenvalues and the sign of the $n+1$ eigenvalue is the same as the sign of $\alpha - v^TA^{-1}v$.
\begin{proof}
Since $A$ is a SPD (SND) matrix then matrix $M$ can be factorized in the following way \cite[Sec.~3.4]{Benzi:2005}
\[
M
=
\begin{bmatrix}
I & 0 \\
v^TA^{-1} & 1
\end{bmatrix}
\begin{bmatrix}
A & 0 \\
0 & \alpha - v^TA^{-1}v
\end{bmatrix}
\begin{bmatrix}
I & A^{-1}v \\
0 & 1
\end{bmatrix}\,.
\]
Here, we have $M = QDQ^T$ and the inertia of the matrix $D$ is the same as the inertia of the matrix $M$.
\end{proof}
\end{lemma}
Lemmas \ref{lem:blockHlin_eig} and \ref{lem:H1HN_lin} tell us that for linear ODEs the Hessian matrix is singular with both positive and negative eigenvalues. Moreover, the dimension of the nullspace of $H$ is at least $(N-2)(n-1)$.

As discussed in \cite{EstrinGreif:2015} higher nullity of $H$ than $(N-1)n+2$, assuming $B(\chi)$ has full column rank, implies that the saddle point matrix \eqref{eq:KKTsystem} is singular. Under an additional assumption on $v_i$ and $\alpha_i$ over Theorem \ref{lem:linH} we are able to conclude that the maximum nullity of $H$ is less then $(N-1)n+2$.
\begin{theorem}
\label{lem:bdHess_lin}
Let $H \in \RR^{N(n+1)\times N(n+1)}$ be a block diagonal matrix with blocks $H_i \in \RR^{(n+1)\times (n+1)}$ $1 \leq i \leq N$. Let each block be of the form
\[
H_i
=
\tbts{A_i}{v_i}{\alpha_i}\,.
\]
Let $A_1 \in \RR^{n \times n}$ and $A_N \in \RR^{n \times n}$ be each either SPD or SND matrix. Let $A_i \in \RR^{n \times n}$, $2 \leq i \leq N-1$, be zero matrices. Moreover, assume that at least one of $v_i$ and $\alpha_i$, $2 \leq i \leq N-1$, in each block is nonzero. Then the maximum dimension of the null-space of $H$ is $(N-2)n + 2$.
\begin{proof}
One can count the possible maximum number of zero eigenvalues of $H$. Using Lemma \ref{lem:H1HN_lin} it follows that there are at most one zero eigenvalue in each block $H_1$ and $H_N$. In the remaining blocks $H_i$, $2 \leq i \leq N-1$, one gets the most zero eigenvalues when $v_i = 0$, $2 \leq i \leq N-1$. It follows from Lemma \ref{lem:blockHlin_eig} that we get additional $(N-2)n$ zero eigenvalues.
\end{proof}
\end{theorem}
For nonlinear ODEs we do not have similar results to those in Lemmas~\ref{lem:blockHlin_eig},~\ref{lem:H1HN_lin} and Theorem~\ref{lem:bdHess_lin}. In our experiments when we used numerical differentiation to compute the Hessian $H$ it always happened to be indefinite.
\section{Solving the Nonlinear Program}
\label{sec:StNLP}
There are two basic approaches to the computation of solution vector $d$ of the saddle-point system~\eqref{eq:KKTsystem} and  the step length $\alpha > 0$: line-search methods and trust-region methods \cite{Nocedal:2006}.
\subsection{Line-search Methods}
\label{subsec:LS}
Line-search methods require either the upper-left block $H$ of \eqref{eq:KKTsystem} to be SPD or the projection of the Hessian onto the null-space of $B^T(\chi)$ to be an SPD matrix \cite[Sec.~18.4]{Nocedal:2006}. When this requirement is fulfilled then the solution $d_\chi$ is a descent direction. 

One method for the approximation of the second derivatives by an SPD matrix is the \emph{BFGS} method~\cite{Nocedal:2006}. Convergence properties of variable metric matrices generated by BFGS were studied in \cite{Ren-pu:1983}. In our problem we try to solve the saddle-point system that features singular and indefinite Hessian $H$. To our knowledge, there are no results showing that the sequence of variable metric matrices generated by BFGS converges when $H$ is not an SPD matrix. 

In our numerical experiments we observed that BFGS produces SPD approximation of the Hessian, however, these were ill-conditioned and near singular. To this end we should not use methods based on the Schur complement reduction~\cite{Benzi:2005} to solve the saddle-point system~\eqref{eq:KKTsystem}.

Since we consider only equality constraints, we solve the saddle-point system~\eqref{eq:KKTsystem} by the \emph{projected preconditioned} conjugate gradient method~\cite[Alg.~NPCG]{LuksanVlcek:2001}. It features an indefinite (constraint) preconditioner~\cite{LuksanVlcek:2001} for which we set the $(1,1)$ block to be the identity matrix~\cite[p.~81]{Benzi:2005}.

It may happen that $d_\chi$ is not acceptable because it increases the violation of the vector of constraints~\cite[Sec.~15.4]{Nocedal:2006}. When this happens one can reset the Hessian to a diagonal SPD matrix and compute a new $d_\chi$. In Section~\ref{sec:CE} we provide a rule for measuring a sufficient decrease in the value of the objective function and satisfaction of constraints. 

The suitable step size $\alpha > 0$ can then be computed using a merit function. In our implementation we use the merit function from the paper \cite{LuksanVlcek:2001} that is
\begin{equation}
\label{eq:merit_fc}
P(\alpha) = F(\chi + \alpha d_\chi) + (\lambda + d_\lambda)^Tc(\chi + \alpha d_\chi) + \frac{\sigma}{2}\|c(\chi + \alpha d_\chi)\|^2,
\end{equation}
where $\sigma \geq 0$. The value $\sigma = 1$ is suitable for many problems~\cite[p.~279]{LuksanVlcek:2001} and we use is in our implementation. For the stepsize selection and termination we used \emph{Backtracking Line Search} \cite[Alg.~3.1]{Nocedal:2006}. In our implementation the stepsize selection is terminated whenever $P(\alpha) - P(0) \leq \delta \alpha P'(0)$ for $\delta = 10^{-4}$. 
%
%
\subsection{Trust-region Methods}
\label{subsec:TR}
An alternative approach to line-search methods are trust-region methods. Instead of computing the direction $d_\chi$ first and then adjusting the step size $\alpha > 0$, trust region methods choose the radius $\Delta > 0$ first and then the suitable direction $d_\chi$ \cite[Ch.~4]{Nocedal:2006}. Trust-region methods are attractive because they do not require the Hessian $H$ in~\eqref{eq:KKTsystem} to be an SPD matrix~\cite{Nocedal:2006}.

Here, we follow the Byrd-Omojokun approach as described in~\cite[Alg.~18.4]{Nocedal:2006}. First we solve the \emph{vertical} subproblem by the dog-leg method~\cite[Ch.~4]{Nocedal:2006}, and then the \emph{horizontal} subproblem by the Steihaug-Toint conjugate gradient method~\cite{Steihaug:1983}, \cite[Alg.~7.2]{Nocedal:2006}. For the computation of the trust-region radius $\Delta > 0$ we use~\cite[Alg.~4.1]{Nocedal:2006}. In addition we use the same merit function $P(\alpha)$ as in the line-search method.

Since the Hessian $H$ may be indefinite and singular, in addition, equality constraints are considered, we use the following \emph{dog-leg} method~\cite[p.~14]{LuksanVlcekTR:2000} in Alg.~\ref{alg:dogleg}.
\begin{algorithm}
\caption{Dog-leg method (vertical step)}
\label{alg:dogleg}
\begin{algorithmic}
\State INPUT: Jacobian $B$, constraints $c$ and $\Delta > 0$
\State $\vartheta \leftarrow 4/5$
\State $d^C \leftarrow - \left( \|Bc\|^2/\|B^TBc\|^2 \right)Bc$\Comment{Cauchy step}
\State $d^N \leftarrow - B\left( B^TB \right)^{-1}c$\Comment{Newton step}
\If{$\|d^C\| \geq \vartheta\Delta$}
\State $d \leftarrow \left( \vartheta\Delta /\|d^C\| \right)d^C$ 
\ElsIf{$\|d^C\| < \vartheta\Delta < \|d^N\|$}
\State  $\alpha \leftarrow -(d^C)^T(d^N-d^C) + \sqrt{\left( (d^C)^T(d^N-d^C) \right)^2 + (\vartheta\Delta)^2 - \|d^C\|^2}$
\State $d \leftarrow d^C + \alpha(d^N-d^C)$
\Else 
\State $d \leftarrow d^N$ 
\EndIf
\State OUTPUT: the direction $d$
\end{algorithmic}
\end{algorithm}
\subsection{Implementation Details for NPCG}
\label{subsec:FurRem}
\begin{figure}
\centering
\includegraphics[width=0.45\linewidth]{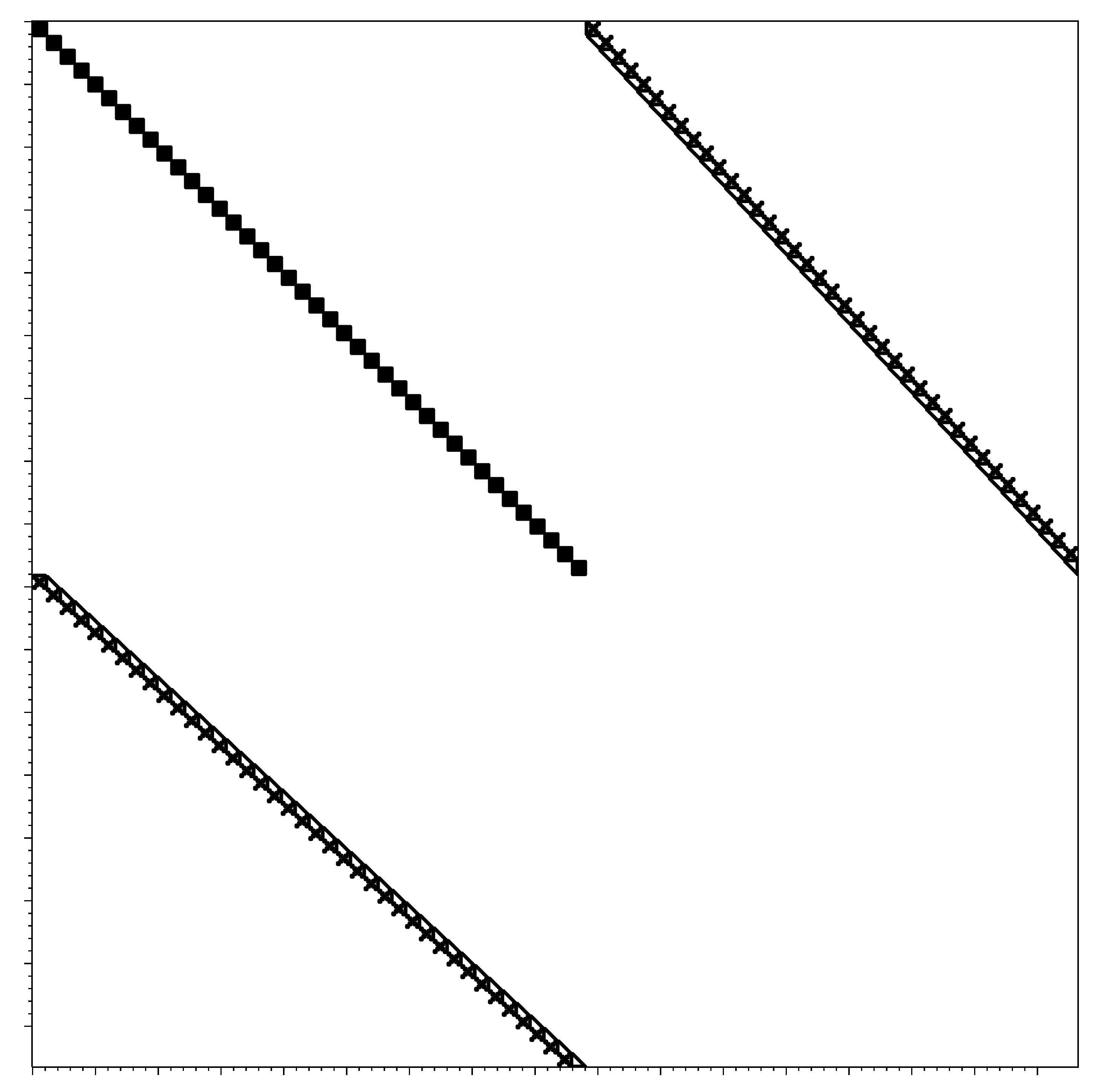}
\includegraphics[width=0.445\linewidth]{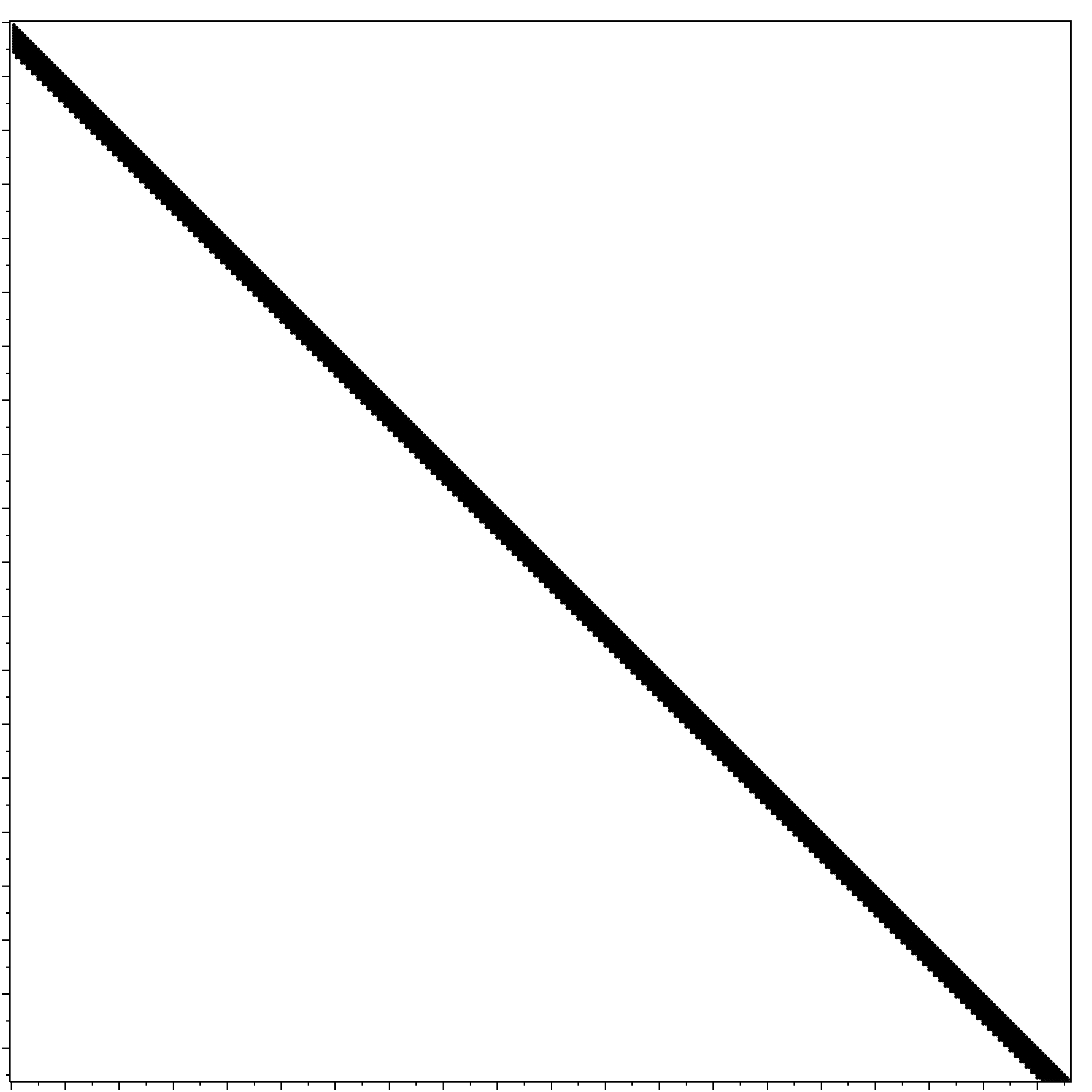}
\caption{The structure of nonzero entries: On the left hand side is the saddle-point matrix $K \in \RR^{832 \times 832}$ with $9250$ nonzero elements. On the right hand side is the Cholesky factor $L \in \RR^{392 \times 392}$ such that $B^TB = LL^T$, with $4447$ nonzero elements. We set $n = 10$, $N = 40$ for the benchmark problem~\ref{subsec:B2} to get these matrices.}
\label{fig:KKT_NNZ}
\end{figure}

Note that the larger the $n$ and $N$ get, the higher the order of the saddle-point matrix $K$ is. However, it is sparse with the structure of nonzero entries as shown in Fig.~\ref{fig:KKT_NNZ}. 
This structure of nonzero elements can be work with efficiently. 

In the projected preconditioned conjugate gradient method~\cite[Alg.~NPCG]{LuksanVlcek:2001} one needs to compute the matrix vector product $Hv$ in each iteration. Since the matrix $H$ is block-diagonal, one should only keep nonzero blocks in memory. Furthermore, since we use the indefinite preconditioner of the form
\[
C
=
\tbts{I}{B}{0}
\]
we need to be able to solve the linear system $B^TBy = z$  efficiently. This is the case as Lemma~\ref{lem:BTB} shows.
\begin{lemma}
\label{lem:BTB}
Let $B$ be the Jacobian of the vector of constraints~\eqref{eq:constrEQ}. Denote by
\begin{align*}
v_\mathrm{I} & =  \elli(x_0^1 - \ceni), \\
 v_\mathrm{U} &= S(t_N, x_0^N)^T\ellu(\Phi(t_N, x_0^N) - \cenu), \\
 \alpha & = \mder{\Phi(t_{N}, x_0^{N})}{t_{N}}^T\ellu(\Phi(t_N, x_0^N) - \cenu),
\end{align*} 
and for $1 \leq i \leq N-1$ put
\begin{align*}
 M_i & = -S(t_i, x_0^i)^T,\\
 v_i &= -\mder{\Phi(t_i, x_0^i)}{t_i},\\
 D_i & = M_i^TM_i + v_iv_i^T + I\,.
\end{align*}
Then 
\begin{equation}
\label{eq:BTB}
B^TB
=
\begin{bmatrix}
v_\mathrm{I}^Tv_\mathrm{I} & v_\mathrm{I}^TM_1 &  &  &  &  & & &  \\ 
M_1^Tv_\mathrm{I} & D_1 & M_2 &  &  &  &  &	& \\  
 & M_2^T &D_2  &  &  &  & 	& \\ 
    &  &  & \ddots &  &  & &	\\ 
   &  &  &  & D_{N-2} & M_{N-1}  &  &	\\ 
 & &  &  & M_{N-1}^T &  D_{N-1} & v_\mathrm{U}	 \\ 
  &  &  &  &  & v_\mathrm{U}^T & v_\mathrm{U}^Tv_\mathrm{U} + \alpha^2
\end{bmatrix}\,.
\end{equation}
\begin{proof}
The formula for $B^TB$ is the result of the matrix multiplication.
\end{proof}
\end{lemma}
The matrix $B^TB$ is SPD and banded. Therefore, the Cholesky factor of $B^TB$ is also not dense and has the banded structure, see Fig.~\ref{fig:KKT_NNZ}. The width of the band in \eqref{eq:BTB} is independent of the number of solution segments $N$. 

Therefore, we do not need to avoid large values of $N$ and could use the benefits of the multiple-shooting methods, hence, restricting the size of intervals of integration~\cite{Ascher:1998}. 
\section{Computational Experiments}
\label{sec:CE}
For computing and approximating the Hessian $H$ we tried different possibilities such as \emph{BFGS}, \emph{SR}-$1$ and  second derivatives computed by finite differences. 
When the true Hessian or its approximation is not an SPD matrix, then line-search may produce a $d_\chi$ that is not a descent direction. In our implementation, we reject the solution vector $d_\chi$  whenever $-P'(0) < 10^{-5} \| d_\chi \| \| \nabla_\chi \LL \|$~\cite[Alg.~3.1]{LuksanVlcek:2001}. In that case we restart the method setting the Hessian $H$ to be the identity matrix. We also set the initial approximation of the Hessian to the identity matrix at the beginning of the optimization.


Our methodology is the following. For the given dynamical system \eqref{eq:DiffEq} and state space dimension $n$ we set $\ceni = [1, \ldots, 1]^T \in \RR^n$ to be the initial state. Then we compute $\cenu = \Phi(5, \ceni)$. The sets $\init$ and $\unsafe$ are balls of radius $1/4$ centred at $\ceni$ and $\cenu$ respectively. We create $N$ solution segments by splitting the solution segment from $\ceni$ to $\cenu$ such that each segment has length $t = 5/N$. Denote the initial states by $x_0^i$, $1 \leq i \leq N$, and modify them by $x_0^i + u$, where  $u = 0.5\times\left[-1, 1, \ldots, (-1)^n\right]^T \in \RR^n$. With these updated initial conditions and lengths $t_i = 5/N$ we get a vector of parameters \eqref{eq:parameters}. From these $N$ segments we try to compute a solution with $x_0^1 \in \init$ and $\Phi(t_N, x_0^N) \in \unsafe$.

We use the following \emph{stopping criteria}: $\|\nabla_\chi \LL \| < 10^{-3}$ and $\| \nabla_\lambda \LL \| < 10^{-8}$ both satisfied (S \textbf{1}); maximum number of iterations is $400$ (S \textbf{2});  the step size $\alpha < 10^{-8}$ for line search and the radius $\Delta < 10^{-8}$ for trust-region (S \textbf{3}); 

If any of these stopping criteria are met then the method terminates. In the end we verify by simulation that parameters \eqref{eq:parameters} give the desired solution. The solution vector $\chi$ is said to be verified by simulation if $\|x_0^1 - \ceni \|_{\elli}^2 < 1 + \varepsilon$ and  $\|\Phi(\sum t_i, x_0^1) - \cenu)\|_{\ellu}^2  < 1 + \varepsilon$, where $\varepsilon = 10^{-4}$. In case they do not it is marked by ``\textbf{F}''. 

In tables with results there are rows and columns marked by an ``S''. In those it is written which stopping criterion took place and values range from $1$ to $3$ as explained above.
  Finally, the number of iterations is denoted by \emph{NIT}.

All computations were carried out in \emph{Scilab} $5.5.2$ \cite{Scilab} installed on a computer with \emph{Cent OS 6.8}. For solving differential equations we used the built-in function \emph{ode}, that in default settings calls the \emph{lsoda} solver of package ODEPACK. As for the sensitivity function~\eqref{eq:Sensitivity}, we either solve the variational equation~\cite[Ch.~7]{Ascher:1998} or use the finite difference method according to the internal differentiation principle~\cite[p.~117]{Betts:2010}. Whenever we needed to solve a system of linear equations, for example obtaining the Newton step in Alg.~\ref{alg:dogleg}, we called the \emph{backslash} operator in Scilab. In the end we applied the same rules for taking and skipping the updates of BFGS and SR-1 methods as in~\cite[Ch.~6]{Nocedal:2006}.

In our implementation, one iteration in line search SQP takes similar amount of running time as one iteration in trust-region SQP. Therefore, we only list the number of iterations in the tables with results. In the end, we did consider inequalities $c_I(\chi) \leq 0$ and used the interior-point method~\cite{Luksan:2004}, however, we did not obtain better results. In general, the interior-point method required more iterations than line-search SQP.

\subsection{Benchmark 1}
\label{subsec:B1}
Consider the following nonlinear dynamical system
\begin{align*}
\dot{x}_1 & = -x_2 + x_1x_3\,, \\
\dot{x}_2 & = x_1 + x_2x_3\,, \\
\dot{x}_3 & = -x_3 - (x_1^2 + x_2^2) + x_3^2\,,
\end{align*}
that is adopted from  from \cite[p.~334]{Khalil:2002}. 
In addition we compare what approximation scheme for the Hessian may be the most suited. 

In Tab. \ref{tab:B1BFGS} there are results when the Hessian blocks $H_1$, $1 \leq i \leq N$, are approximated by BFGS. When we used SR-$1$ we got results in Tab. \ref{tab:B1SR1}. Note that in cases where the maximum number of iterations was reached we got feasibility conditions $\| \nabla_\lambda \LL \| \leq 10^{-7}$ except for $N = 25$ where  $\| \nabla_\lambda \LL \| \leq 10^{-4}$. For the results in Tab. \ref{tab:B1ExactNonlin} a built-in function \emph{numderivative} was used for computing the Hessian. Both approaches yield similar results in terms of the number of iterations in this case.
\begin{table}
\centering
\begin{tabular}{l|*{6}{c}}
N & $5$  &  $10$  &  $15$  &  $20$  &  $25$  &  $30$ \\
NIT & $35$  &  $33$  &  $31$  &  $48$  &  $47$  &  $49$ \\
S & $1$  &  $1$  &  $1$  &  $1$  &   $3$  &  $3$ \\ 
\end{tabular}
\qquad
\begin{tabular}{l|*{6}{c}}
 N & $5$  &  $10$  &  $15$  &  $20$  &  $25$  &  $30$ \\
 NIT & $28$  &  $32$  &  $36$  &  $43$  &  $45$  &  $51$ \\
 S & $1$  &  $1$  &  $1$  &  $1$  &  $1$  &  $1$ \\
\end{tabular}
\caption{Comparison of line-search and trust region on benchmark~\ref{subsec:B1} where the Hessian $H$ is approximated by the BFGS method. Left: line-search; Right: trust region.}
\label{tab:B1BFGS}
\end{table}
\begin{table}
\centering
\begin{tabular}{l|*{6}{c}}
 N & $5$  &  $10$  &  $15$  &  $20$  &  $25$  &  $30$ \\
 NIT & $16$  &  $28$  &  $34$  &  $46$  &  $7$  & $47$ \\
 S & $1$  &  $1$  &  $3$  &  $3$  &  F  &  $3$ \\
\end{tabular}
\quad
\begin{tabular}{l|*{6}{c}}
 N & $5$  &  $10$  &  $15$  &  $20$  &  $25$  &  $30$ \\
 NIT & $29$  &  $400$  &  $127$  &  $400$  &  $400$  &  $400$ \\
S &  $1$  &  $2$  &  $1$  &  $2$  &  $2$  &  $2$ \\ 
\end{tabular}
\caption{Comparison of line-search and trust region on benchmark~\ref{subsec:B1} where the Hessian $H$ is approximated by the SR1 method. Left: line-search; Right: trust region. In one instance line-search failed at finding a desired solution.}
\label{tab:B1SR1}
\end{table}
\begin{table}
\centering
\begin{tabular}{l|*{6}{c}}
N &  $5$  &  $10$  &  $15$  &  $20$  &  $25$  &  $30$ \\
NIT &  $29$  &  $41$  &  $20$  &  $25$  &  $31$ &  $24$ \\
 S & $1$  &  $1$  &  $3$  &  $3$  &  $3$  &  $3$ \\
\end{tabular}
\qquad
\begin{tabular}{l|*{6}{c}}
 N & $5$  &  $10$  &  $15$  &  $20$  &  $25$  &  $30$ \\
NIT &  $24$  &  $44$  &  $21$  &  $27$  &  $46$ &  $36$ \\
 S & $1$  &  $1$  &  $1$  &  $1$  &  $1$  &  $1$ \\ 
\end{tabular}
\caption{Comparison of line-search and trust region on benchmark~\ref{subsec:B1} where the Hessian $H$ is computed by \emph{numderivative}~\cite{Scilab}. Left: line-search; Right: trust region.}
\label{tab:B1ExactNonlin}
\end{table}

Let us go back to the results in Tab. \ref{tab:B1SR1} where we used SR-$1$ method to approximate the Hessian. We list the number of restarts  in the line-search method when we needed to change the approximated Hessian for the identity matrix to get a descent direction $d_\chi$. With increasing $N$ we have: $1, 0, 2, 4, 0$ and $2$ restarts.
\subsection{Benchmark 2}
\label{subsec:B2}
Consider the following linear dynamical system
\begin{align*}
\dot{x} & = Ax \\
 & =  
\begin{bmatrix}
\tbt{0}{1}{-1}{0} & & \\
& \ddots & \\
& & \tbt{0}{1}{-1}{0}
\end{bmatrix}x\,,
\end{align*}
where $A \in \RR^{n \times n}$. This benchmark problem can not only be scaled up in the number of solution segments $N$ but as well as in the state space dimension $n$. When $n = 40$ and $N = 30$ we solve constrained optimization problem with $N(n+1) = 1230$ parameters. 


In Tab. \ref{tab:B2BFGS} there are results when the BFGS method was used. One can see that the line search method outperforms the trust-region method, especially, for the higher values of $n$. The results for SR-$1$ approximation scheme are shown in Tab.~\ref{tab:B2SR1}. Note that, when SR-$1$ was used, the number of restarts in line search method for $n = 40$ and $N = 5$ was zero. However, for the setting $n = 40$ and $N = 28$ we needed to reset the Hessian to be the identity matrix twenty-eight times.

When the formulas from Theorem~\ref{lem:linH} are used one gets the results in Tab.~\ref{tab:B2ExactLin}. Notice that the trust-region approach almost always terminates because of the maximum number of iterations condition. However, when we investigate the feasibility condition we get that $\| \nabla_\lambda \LL \| \leq 10^{-4}$ in all of those cases. In particular for $n = 40$ and $N = 10$ we have $\|\nabla_\chi \LL \| \leq 10^{-2}$ and $\| \nabla_\lambda \LL \| \leq 10^{-12}$. 
\begin{table}
\centering
\begin{tabular}{*{24}{c}}
n & N & NIT & S & n & N & NIT & S & n & N & NIT & S & n & N & NIT & S \\
\hline
 $10$  &  $5$  &  $28$  &  $1$  &  $20$  &  $5$  &  $33$  &  $1$  &  $30$  &  $5$  &  $30$  & $1$  &  $40$  &  $5$  &  $34$  &  $1$ \\
   &  $10$  &  $31$  &  $1$  &    &  $10$  &  $39$  &  $1$  &    &  $10$  &  $172$  &  $1$  &   &  $10$  &  $29$  &  $1$ \\
   &  $15$  &  $400$  &  $2$  &    &  $15$  &  $36$  &  $1$  &   &  $15$ &  $39$  &  $1$  &    &  $15$  &  $37$  &  $1$ \\
   &  $20$  &  $48$  &  $1$  &   &  $20$  &  $41$  &  $1$  &   &  $20$  &  $172$  &  $1$  &    &  $20$  &  $32$  &  $1$ \\
   &  $25$  &  $42$  &  $1$  &    &  $25$  &  $36$  &  $1$  &    &  $25$  &  $108$  &  $1$  &    &  $25$  &  $87$  &  $3$ \\
   &  $30$  &  $35$  &  $1$  &    &  $30$  &  $39$  &  $1$  &    &  $30$  &  $44$  &  $1$  &    &  $30$  &  $51$  &  $1$ \\[5mm]
n & N & NIT & S & n & N & NIT & S & n & N & NIT & S & n & N & NIT & S \\
\hline
 $10$  &  $5$  &  $26$  &  $1$  &  $20$  &  $5$  &  $28$  &  $1$  &  $30$  &  $5$  &  $29$  &  $1$  &  $40$  &  $5$  &  $30$  &  $1$ \\
&  $10$  &  $91$  &  $1$  &    &  $10$  &  $400$  &  $2$  &    &  $10$  &  $400$  &  $2$  &    &  $10$  &  $400$  &  $2$ \\
&  $15$  &  $119$  &  $1$  &    &  $15$  &  $90$  &  $1$  &    &  $15$  &  $308$  &  $1$  &    &  $15$  &  $400$  &  $2$ \\
&  $20$  &  $400$  &  $2$  &    &  $20$  &  $101$  &  $1$  &    &  $20$  &  $400$  &  $2$  &    &  $20$  &  $400$  &  $2$ \\
&  $25$  &  $38$  &  $1$  &    &  $25$  &  $95$  &  $1$  &    &  $25$  &  $400$  &  $2$  &    &  $25$  &  $400$  &  $2$ \\
&  $30$  &  $39$  &  $1$  &    &  $30$  &  $53$  &  $1$  &    &  $30$  &  $112$  &  $1$  &    &  $30$  &  $311$  &  $1$ \\ 
\end{tabular}
\caption{Comparison of line-search and trust region on benchmark~\ref{subsec:B2} where the Hessian $H$ is approximated by the BFGS method. Top: line-search; Bottom: trust region.}
\label{tab:B2BFGS}
\end{table}
\begin{table}
\centering
\begin{tabular}{*{24}{c}}
n & N & NIT & S & n & N & NIT & S & n & N & NIT & S & n & N & NIT & S \\
\hline
 $10$  &  $5$  &  $64$  &  $1$  &  $20$  &  $5$  &  $35$  &  $1$  &  $30$  &  $5$  &  $39$  & $1$  &  $40$  &  $5$  &  $27$  &  $1$ \\
 &  $10$  &  $37$  &  $1$  &  &  $10$  &  $34$  &  $1$  & &  $10$  &  $35$  &  $1$  & &  $10$  &  $34$  &  $1$ \\
 &  $15$  &  $53$  &  $1$  &  &  $15$  &  $70$  &  $1$  & &  $15$ &  $50$  &  $1$  & &  $15$  &  $62$  &  $1$ \\
 &  $20$  &  $38$  &  $1$  &  &  $20$  &  $34$  &  $1$  & &  $20$  &  $58$  &  $1$  & &  $20$  &  $45$  &  $1$ \\
 &  $25$  &  $52$  &  $1$  &  &  $25$  &  $41$  &  $1$  & &  $25$  &  $45$  &  $1$  & &  $25$  &  $128$  &  $1$ \\
 &  $30$  &  $54$ &  $1$  &  &  $30$  &  $46$  &  $1$  & &  $30$  &  $100$  &  $1$  &&  $30$  &  $57$  &  $1$ \\[5mm]
n & N & NIT & S & n & N & NIT & S & n & N & NIT & S & n & N & NIT & S \\
\hline
 $10$  &  $5$  &  $400$  &  $2$  &  $20$  &  $5$  &  $400$  &  $2$  &  $30$  &  $5$  &  $400$ &  F  &  $40$  &  $5$  &  $400$  &  F \\
  &  $10$  &  $400$  &  $2$  &    & $10$  &  $400$  &  $2$  &    &  $10$  &  $400$  &  $2$  &    &  $10$  &  $400$ &  $2$ \\
  &  $15$  &  $400$  &  $2$  &    &  $15$  &  $400$  &  $2$  &   &  $15$  &  $400$  &  $2$  &    &  $15$  &  $400$  &  $2$ \\
  &  $20$  &  $400$ &  F       &    &  $20$  &  $400$  &  $2$  &    &  $20$  &  $400$  &  F      &   &  $20$  &  $400$  &  F \\
  &  $25$  &  $400$  &  F      &    &  $25$  &  $400$  &  $2$  &    &  $25$  &  $400$  &  $2$  &    &  $25$  &  $400$  &  $2$ \\
  &  $30$  &  $400$  &  $2$  &    &  $30$  &  $400$  &  $2$  &    &  $30$  &  $400$  &  $2$  &    &  $30$  &  $400$  &  $2$ \\
\end{tabular}
\caption{Comparison of line-search and trust region on benchmark~\ref{subsec:B2} where the Hessian $H$ is approximated by the SR1 method. Top: line-search; Bottom: trust region. There are six failed attempts at finding a desired solution for the trust-region method.}
\label{tab:B2SR1}
\end{table}
\begin{table}
\centering
\begin{tabular}{*{24}{c}}
n & N & NIT & S & n & N & NIT & S & n & N & NIT & S & n & N & NIT & S \\
\hline
 $10$  &  $5$  &  $43$  &  $1$  &  $20$  &  $5$  &  $47$  &  $1$  &  $30$  &  $5$  &  $27$  & $1$  &  $40$  &  $5$  &  $26$  &  $1$ \\
   &  $10$  &  $117$  &  $1$  &   &  $10$  &  $108$  &  $1$  &    &  $10$  &  $121$  &  $1$  &    &  $10$  &  $23$  &  $1$ \\
   &  $15$  &  $49$  &  $1$  & &  $15$  &  $32$  &  $1$  &   &  $15$  &  $76$  &  $1$  &   &  $15$  &  $41$  &  $1$ \\
   &  $20$  &  $47$  &  $1$  & &  $20$  &  $23$  &  $3$  &   &  $20$  &  $48$  &  $1$  &   &  $20$ &  $30$  &  $1$ \\
   &  $25$  &  $29$  &  $1$  & &  $25$  &  $49$  &  $1$  &  &  $25$  &  $39$  &  $1$  &  &  $25$  &  $31$  &  $1$ \\
   &  $30$  &  $32$  &  $1$  & &  $30$  &  $34$  &  $1$  &   &  $30$  &  $34$  &  $1$  &   &  $30$  &  $33$  &  $1$ \\[5mm]
n & N & NIT & S & n & N & NIT & S & n & N & NIT & S & n & N & NIT & S \\
\hline
 $10$  &  $5$  &  $185$  &  $1$  &  $20$  &  $5$  &  $28$  &  $1$  &  $30$  &  $5$  &  $400$ &  $2$  &  $40$  &  $5$  &  $400$  &  $2$ \\
 &  $10$  &  $400$  &  $2$  &    & $10$  &  $400$  &  $2$  &    &  $10$  &  $400$  &  $2$  &    &  $10$  &  $400$ &  $2$ \\
 &  $15$  &  $247$  &  $1$  &   &  $15$  &  $400$  &  $2$  &    &  $15$  &  $400$  &  $2$  &    &  $15$  &  $400$  &  $2$ \\
 &  $20$  &  $400$  &  F      &    &  $20$  &  $400$  &  $2$  &    &  $20$  &  $400$  &  $2$  &   &  $20$  &  $400$  &  F \\
 &  $25$  &  $400$  &  $2$  &    &  $25$  &  $400$ &  $2$  &   &  $25$  &  $400$  &  $2$  &    &  $25$  &  $400$  &  $2$ \\
 &  $30$  &  $400$  &  $2$  &   &  $30$  &  $400$  &  $2$  &  &  $30$  &  $400$  &  $2$  &    &  $30$  &  $400$  &  $2$ \\ 
\end{tabular}
\caption{Comparison of line-search and trust region on benchmark~\ref{subsec:B2} where the Hessian $H$ is computed using formulas from Theorem~\ref{lem:linH}. Top: line-search; Bottom: trust region. There are two failed attempts at finding a desired solution for the trust-region method.}
\label{tab:B2ExactLin}
\end{table}
%
%
%
%
%
%
\subsection{Benchmark 3}
\label{subsec:B3}
Consider the following nonlinear dynamical system
\begin{align*}
\dot{x} & = Ax + \sin(x^r) \\
 & =  
\begin{bmatrix}
\tbt{0}{1}{-1}{0} & & \\
& \ddots & \\
& & \tbt{0}{1}{-1}{0}
\end{bmatrix}x
+
\begin{bmatrix}
\sin (x_n) \\
\vdots \\
\sin (x_1)
\end{bmatrix}\,,
\end{align*}
where $A \in \RR^n$. It is similar to benchmark \ref{subsec:B2}, however, this time there is a nonlinear term $\sin(x^r)$ present. This causes that blocks $A_i$, $1 \leq i \leq N$, in the Hessian \ref{lem:nonlinH} to be nonzero in general.

One can see in Tab.~\ref{tab:B3BFGS} that both approaches yield similar results when BFGS was used. All runs terminated successfully with no fails. As for the SR-$1$ approximation of the Hessian $H$ the results in Tab.~\ref{tab:B3SR-1} are not that promising. Both approaches failed to find an acceptable solution in few cases.

Note that in Tab. \ref{tab:B3SR-1} the trust-region method always used the maximum number of iterations. The reason behind is that the norm $\| \nabla_\chi \LL \|$ does not drop below the prescribed tolerance $10^{-3}$. 
In addition, from our experience it follows that the number of restarts in line search does not tell us whether our method converges to a desired solution. For instance, when SR-$1$ was used, there was only one restart for $n = 10$ and $N = 25$, yet our method failed to find a solution. On the other hand our method found a solution for $n = 40$ and $N = 25$ although it had to reset the Hessian for the identity matrix eighteen times.
%
%
%
%
%
%
\begin{table}
\centering
\begin{tabular}{*{24}{c}}
 n & N & NIT & S & n & N & NIT & S & n & N & NIT & S & n & N & NIT & S \\
\hline
 $10$  &  $5$  &  $37$  &  $1$  &  $20$  &  $5$  &  $54$  &  $1$  &  $30$  &  $5$  &  $46$  & $1$  &  $40$  &  $5$  &  $67$  &  $1$ \\
  &  $10$  &  $47$  &  $1$  &   &  $10$  &  $54$  &  $3$  &   &  $10$  &  $122$  &  $1$  &  &  $10$  &  $57$  &  $1$ \\
  &  $15$  &  $126$  &  $1$  &   &  $15$  &  $52$  &  $1$  &   &  $15$ &  $97$  &  $1$  &  &  $15$  &  $60$  &  $1$ \\
  &  $20$  &  $79$  &  $3$  & &  $20$  &  $149$  &  $3$  &   &  $20$  &  $101$  &  $3$  &  &  $20$ &  $90$  &  $3$ \\
  &  $25$  &  $185$  &  $1$  &   &  $25$  &  $107$  &  $3$ &   &  $25$  &  $80$  &  $1$  &  &  $25$  &  $98$  &  $3$ \\
  &  $30$  & $400$  &  $2$  &   &  $30$  &  $108$  &  $1$  &   &  $30$  &  $122$  &  $1$  &  &  $30$  &  $98$  &  $1$ \\[5mm]
n & N & NIT & S & n & N & NIT & S & n & N & NIT & S & n & N & NIT & S \\
\hline
 $10$  &  $5$  &  $80$  &  $1$  &  $20$  &  $5$  &  $110$  &  $1$  &  $30$  &  $5$  &  $81$  & $1$  &  $40$  &  $5$  &  $57$  &  $1$ \\
  &  $10$  &  $54$  &  $1$  & &  $10$  &  $54$  &  $1$  &    &  $10$  &  $149$  &  $1$  &    &  $10$  &  $118$  &  $1$ \\
  &  $15$  &  $152$  &  $1$  & &  $15$  &  $126$  &  $1$  &    &  $15$  &  $111$  &  $1$  &    &  $15$  &  $56$  &  $1$ \\
  &  $20$  &  $121$  &  $1$  & &  $20$  &  $133$  &  $1$  &    &  $20$  &  $89$  &  $1$  &    &  $20$  &  $192$  &  $1$ \\
  &  $25$  &  $61$  &  $1$  & &  $25$  &  $63$  &  $1$  &    &  $25$  &  $69$  &  $1$  &    &  $25$  &  $165$  &  $1$ \\
  &  $30$ &  $152$  &  $1$  & &  $30$  &  $68$  &  $1$  &    &  $30$  &  $178$  &  $1$  &   &  $30$  &  $142$  &  $1$ \\ 
\end{tabular}
\caption{Comparison of line-search and trust region on benchmark~\ref{subsec:B3} where the Hessian $H$ is approximated by the BFGS method. Top: line-search; Bottom: trust region.}
\label{tab:B3BFGS}
\end{table}
\begin{table}
\centering
\begin{tabular}{*{24}{c}}
 n & N & NIT & S & n & N & NIT & S & n & N & NIT & S & n & N & NIT & S \\
\hline
 $10$  &  $5$  &  $65$  &  $1$  &  $20$  &  $5$  &  $31$  &  $1$  &  $30$  &  $5$  &  $49$  & $1$  &  $40$  &  $5$  &  $52$  &  $1$ \\
  &  $10$  &  $68$  &  $1$  &   &  $10$  &  F  &  $2$  &    &  $10$  &  $37$  &  $3$  &    &  $10$  &  $51$  &  $1$ \\
  &  $15$  &  $40$  &  $3$  &   &  $15$  &  $75$  &  $1$  &  &  $15$ &  $60$  &  $1$  &   &  $15$  &  $54$  &  $1$ \\
  &  $20$  &  $90$  &  F       &  &  $20$  &  $75$  &  $3$  &  &  $20$  &  $82$  &  $1$  &   &  $20$  & $92$  &  $3$ \\
  &  $25$  &  $21$  &  F       &   &  $25$  &  F  &  $3$  &  &  $25$  &  $101$  &  $3$  &   &  $25$  &  $99$  &  $1$ \\
  &  $30$  &  $88$  &  $1$  &   &  $30$  &  $86$  &  $1$  &  &  $30$  &  $122$  &  $1$  &   &  $30$  &  $80$  &  $1$ \\ [5mm]
n & N & NIT & S & n & N & NIT & S & n & N & NIT & S & n & N & NIT & S \\
\hline
 $10$  &  $5$  &  -      &  F     &  $20$  &  $5$  &  $400$  &  $2$  &  $30$  &  $5$  &  $400$  &  F     &  $40$  &  $5$  &  $400$  &  F   \\
  &  $10$  &  $400$  &  $2$  &  &  $10$ &  $400$  &  $2$  &   &  $10$  &  $400$  &  $2$  & &  $10$  &  $400$  &  $2$ \\
  &  $15$  &  $400$  &  $2$  &  &  $15$  &  $400$  &  $2$  &   &  $15$  &  $400$  &  $2$  & &  $15$  &  $400$  &  $2$ \\
  &  $20$  &  $400$  &  $2$  &  &  $20$  &  $400$  &  $2$  &   &  $20$  &  $400$  &  $2$  & & $20$  &  $400$  &  $2$ \\
  &  $25$  &  $400$  &  $2$  &  &  $25$  &  $400$  &  F      &  & $25$  &  $400$  &  $2$  &    &  $25$  &  $400$  &  $2$ \\
  & $30$  &  $400$  &  $2$  &  &  $30$  &  $400$  &  $2$  &   &  $30$  &  $400$ &  $2$  &   &  $30$  &  $400$  &  F \\ 
\end{tabular}
\caption{Comparison of line-search and trust region on benchmark~\ref{subsec:B3} where the Hessian $H$ is approximated by the SR1 method. Top: line-search; Bottom: trust region. There are failed attempts at finding a desired solution for both methods.}
\label{tab:B3SR-1}
\end{table}
\section{Conclusion}
\label{sec:End}
In this paper we investigated the problem of finding a solution to a dynamical system that starts in a given set and reaches another set of states. We studied properties of the saddle-point matrix \eqref{eq:KKTsystem} resulting from the minimization formulation~\eqref{eq:MinProbEQ}. In addition, we compared line search and trust-region methods on benchmark problems. We conclude that the most suitable approach to solving the minimization problem~\eqref{eq:MinProbEQ} is line-search SQP. In more detail the most promising approach to solve such problem
\begin{itemize}
\item uses BFGS block-wise to approximate $\nabla_\chi^2 \LL$, and
\item solves the saddle-point system by projected preconditioned CG.
\end{itemize}

As for the properties of the saddle-point matrix Theorems~\ref{lem:nonlinH} and~\ref{lem:linH} show the structure of nonzero elements of the Hessian $\nabla_\chi^2 \LL$. Moreover, for the linear dynamic one can use the formulas from Theorem \ref{lem:linH} to compute second derivatives fast using already computed data from the Jacobian of constraints. To a lesser extent the same can be said in the nonlinear case about the result in Theorem~\ref{lem:nonlinH}, where formulas for some elements in the Hessian are provided. In general, in this problem the upper-left block of the saddle-point matrix $K$ is singular and indefinite.
%
%
%
\bibliographystyle{abbrv}
\bibliography{../bibliography/kuratko}

\end{document}